\documentclass[12pt]{article}

\usepackage{amssymb,latexsym}

\newtheorem{theorem}{Theorem}[section]
\newtheorem{proposition}[theorem]{Proposition}
\newtheorem{corollary}[theorem]{Corollary}
\newtheorem{definition}[theorem]{Definition}

\newenvironment{proof}{\medskip\noindent{\it Proof.\ }}{\mbox{$\Box$}\medskip}

\begin{document}

\def\cQ{{\cal Q}}
\def\II{{\mathfrak A}}
\def\JJ{{\mathfrak C}}
\def\vr{\varnothing}

\title{231-Avoiding Involutions and Fibonacci Numbers}

\author{Eric S. Egge \\
Department of Mathematics \\
Gettysburg College\\
Gettysburg, PA  17325  USA \\[4pt]
eggee@member.ams.org \\
\\
Toufik Mansour\\
Department of Mathematics\\
Chalmers University of Technology\\
412 96 G\"oteborg\\
Sweden
\\[4pt]
toufik@math.chalmers.se}
\maketitle

\begin{abstract}
We use combinatorial and generating function techniques to enumerate various sets of involutions which avoid 231 or contain 231 exactly once.
Interestingly, many of these enumerations can be given in terms of $k$-generalized Fibonacci numbers.

\medskip

\end{abstract}

\section{Introduction and Notation}

Let $S_n$ denote the set of permutations of $\{1, \ldots, n\}$, written in one-line notation, and suppose $\pi \in S_n$.
We say $\pi$ is an {\it involution} whenever $\pi(\pi(i)) = i$ for all $i$, $1 \le i \le n$, and we write $I_n$ to denote the set of involutions in $S_n$.
Now suppose $\pi \in S_n$ and $\sigma \in S_k$.
We say $\pi$ {\it avoids} $\sigma$ whenever $\pi$ contains no subsequence with all of the same pairwise comparisons as $\sigma$.
For example, the permutation 214538769 avoids 312 and 2413, but it has 2586 as a subsequence so it does not avoid 1243.
If $\pi$ avoids $\sigma$ then $\sigma$ is sometimes called a {\it pattern} or a {\it forbidden subsequence} and $\pi$ is sometimes called a {\it restricted permutation} or a {\it pattern-avoiding permutation}.
In this paper we will be interested in permutations which avoid several patterns, so for any set $R$ of permutations we write $S_n(R)$ to denote the elements of $S_n$ which avoid every element of $R$.
For any set $R$ of permutations we take $S_n(R)$ to be the empty set whenever $n < 0$ and we take $S_0(R)$ to be the set containing only the empty permutation.
When $R = \{\pi_1, \pi_2, \ldots, \pi_r\}$ we often write $S_n(R) = S_n(\pi_1, \pi_2, \dots, \pi_r)$.
We will also be interested in involutions which avoid several patterns, so for any set $R$ of permutations we write $I_n(R)$ to denote the set of involutions in $S_n(R)$.

In this paper we will encounter several sequences, many of which can be written in terms of a particular family of sequences.
To define this family of sequences, first fix $k \ge 0$.
For all $n \le 0$ we set $F_{k,n} = 0$;  we also set $F_{k,1} = 1$.
For all $n \ge 2$ we set $F_{k,n} = \sum_{i=1}^k F_{k,n-i}$.
We observe that $F_{2,n}$ is the ordinary Fibonacci number $F_n$ for all $n \ge 0$.
In view of this, we refer to the numbers $F_{k,n}$ as the {\it $k$-generalized Fibonacci numbers}.
We also observe that the ordinary generating function for the $k$-generalized Fibonacci numbers is given by
\begin{equation}
\label{eqn:kfibgf}
\sum_{n=0}^\infty F_{k,n} x^n = \frac{x}{1 - x - x^2 - \ldots - x^k} \hspace{30pt} (k \ge 1).
\end{equation}
We will also make use of the fact, easily verified by induction, that the number of tilings of a $1 \times n$ rectangle with tiles of size $1 \times 1$, $1 \times 2, \ldots, 1 \times k$ is given by $F_{k,n+1}$ for all $k \ge 1$ and all $n \ge 0$.

Restricted permutations were first connected with Fibonacci numbers by Simion and Schmidt, who showed in \cite[Prop. 15]{SimionSchmidt} that
$$|S_n(132, 213, 123)| = F_{n+1} \hspace{30pt} (n \ge 0).$$
The present authors generalized this result extensively in \cite{EggeMansourFib}, where they gave several families of sets of restricted permutations which can be counted in terms of $k$-generalized Fibonacci numbers.
For example, the present authors showed in \cite{EggeMansourFib} that
\begin{equation}
\label{eqn:EM1}
|S_n(132, 213, \beta_{a,b,c})| = \sum_{k=1}^{a+c-1} {{n-1} \choose {k-1}} + \sum_{k=a+c}^n {{k-1} \choose {a+c-1}} F_{b-1,n-k+1},
\end{equation}
where $\beta_{a,b,c}$ is the permutation in $S_{a+b+c}$ given by
$$\beta_{a,b,c} = a+b+c, a+b+c-1, \ldots, b+c+1, c+1, c+2, \ldots, b+c, c, c-1, \ldots, 2,1.$$
In fact, $|S_n(132, 213, \tau)|$ can be expressed in terms of $k$-generalized Fibonacci numbers for every $\tau \in S_n(132, 213)$, as Mansour demonstrated in \cite{Mansour33k} by expressing the generating function for $|S_n(132, 213, \tau)|$ as a determinant of a matrix of generating functions for various $k$-generalized Fibonacci numbers.

Although they were initially studied through the Robinson-Schensted correspondence, restricted involutions have recently begun to receive attention as objects of study in their own right.
In \cite{Regev} Regev provided an asymptotic formula for $|I_n(12\ldots k)|$ and showed that
$|I_n(1234)| = M_n$, where $M_n$ is the $n$th Motzkin number, which may be defined by $M_0 = 1$ and $M_n = M_{n-1} + \sum_{i=0}^{n-2} M_i M_{n-i-2}$ for $n \ge 1$.
In \cite{Gessel} Gessel enumerated $I_n(12\ldots k)$, and Gouyou--Beauchamps \cite{GouyouBYoung} has given bijective proofs of exact formulas for $|I_n(12345)|$ and $|I_n(123456)|$.
Guibert \cite{G} has established bijections between 1-2 trees
with $n$ edges and several sets of restricted involutions, including $I_n(3412)$ and $I_n(4321)$.
This latter bijection leads, by way of the Robinson-Schensted correspondence, to a bijection between 1-2 trees with $n$ edges and $I_n(1234)$.
Guibert has also given \cite{G} a
bijection between vexillary involutions of length $n$ (that is, $I_n(2143)$) and $I_n(1243)$. 
More recently, Guibert, Pergola and Pinzani \cite{GPP} gave a bijection between 1-2 trees with $n$ edges and vexillary involutions of length $n$. 
Combining all of these results, we find that $|I_n(1234)| = |I_n(3412)| = |I_n(4321)| = |I_n(2143)| = |I_n(1243)| = M_n$.
At this writing it is an open problem to prove the conjecture of Guibert \cite{G} that $|I_n(1432)| = M_n$. 

In this paper we use combinatorial and generating function techniques to enumerate various sets of involutions which avoid 231 or contain 231 exactly once.
It turns out that many of these enumerations can be given in terms of $k$-generalized Fibonacci numbers.
In section 2 we use results of Simion and Schmidt \cite{SimionSchmidt}, Mansour \cite{Mansour33k}, and the current authors \cite{EggeMansourFib} to explain how to enumerate involutions which avoid 231 and another pattern.
In section 3 we enumerate involutions which avoid 231 and contain a given pattern.
In section 4 we enumerate involutions which contain 231 exactly once and avoid another pattern.
In section 5 we enumerate involutions which contain 231 exactly once and contain another pattern.

\section{Involutions Which Avoid 231 and Avoid Another Pattern}
\label{sec:In231tau}

In this section we briefly consider $I_n(231, \tau)$, where $\tau \in S_k$.
Our consideration will be brief because these sets of permutations have already been extensively studied in a slightly different guise.
To describe this different guise, we begin with some notation.

\begin{definition}
\label{defn:layeredpermutation}
Fix $n \ge 1$ and let $l_1, l_2, \ldots, l_m$ denote a sequence such that $l_i \ge 1$ for $1 \le i \le m$ and $\sum_{i=1}^m l_i = n$.
We write $[l_1, l_2, \ldots, l_m]$ to denote the permutation given by
$$[l_1, l_2, \ldots, l_m] = l_1, l_1-1, \ldots, 1, l_2 + l_1, l_2+l_1-1, \ldots, l_1+1, \ldots, n, n-1, \ldots, n-l_m+1.$$
We call a permutation {\em layered} whenever it has the form $[l_1,\ldots,l_m]$ for some sequence $l_1, \ldots, l_m$.
\end{definition}

We remark that layered permutations have also been studied in \cite{BonaLayered} and \cite{MansourVainshtein}.

Associating the layered permutation $[l_1, \ldots, l_m]$ with the tiling $1 \times l_1$, $1 \times l_2$, $1 \times l_3, \ldots, 1 \times l_m$ of a rectangle of size $1 \times n$, we obtain a natural bijection between layered permutations of length $n$ and tilings of a $1 \times n$ rectangle with rectangles of size $1 \times 1, 1 \times 2, 1 \times 3, \ldots$.

Combining the proof of \cite[Prop. 6]{SimionSchmidt}, the remarks following this proof, and the remarks following the proof of \cite[Prop. 12]{SimionSchmidt}, we obtain the following connection between $I_n(231)$ and layered permutations.

\begin{proposition}
\label{prop:In231form}
For all $\pi \in I_n$, the number of subsequences of type 231 in $\pi$ is equal to the number of subsequences of type 312 in $\pi$.
Moreover, for all $n \ge 0$, the sets $I_n(231)$, $I_n(312)$, and $S_n(231, 312)$ are all equal to the set of layered permutations of length $n$.
\end{proposition}

It follows from this result that $I_n(231,\tau) = S_n(231, 312, \tau) = I_n(312, \tau)$ for any permutation $\tau$.
Writing $\tau^r$ to denote the permutation obtained by writing the entries of $\tau$ in reverse order, we find that $S_n(231, 312, \tau) = S_n(132, 213, \tau^r)$ for any permutation $\tau$.
These observations allow us to translate results such as (\ref{eqn:EM1}) into enumerations of $I_n(213,\tau)$ for various permutations $\tau$.

\section{Involutions Which Avoid 231 and Contain Another Pattern}

In this section we consider those involutions in $I_n(231)$ which contain a given pattern $\tau$.
We begin by setting some notation.

\begin{definition}
For all $n \ge 0$, all $r \ge 0$, and all permutations $\tau$, we write $I_n^r(231 | \tau)$ to denote the set of involutions in $I_n(231)$ which contain exactly $r$ subsequences of type $\tau$.
We write $\II_\tau^r(x)$ to denote the generating function given by
$$\II_\tau^r(x) = \sum_{n=0}^\infty |I_n^r(231 | \tau)| x^n.$$
We write $\II_\tau(x,y)$ to denote the generating function given by
$$\II_\tau(x,y) = \sum_{n=0}^\infty \sum_{r=0}^\infty |I_n^r(231| \tau)| x^n y^r.$$
\end{definition}

In view of Proposition \ref{prop:In231form}, if $\tau$ is not layered then $\II_\tau(x,y) = 0$, so it is reasonable to ask for a closed form expression for $\II_\tau(x,y)$ for any layered permutation $\tau$. 
Such an expression appears to be difficult to obtain, so we content ourselves here with a closed form expression for $\II_{k\ldots 21}(x,y)$.

\begin{theorem}
\label{thm:avoid231contain12kgf}
For all $k\geq 1$, we have
\begin{equation}
\label{eqn:IIk21gf}
\II_{k\dots 21}(x,y)=\frac{1}{1-\sum\limits_{j\geq 1} x^j y^{{{j}\choose{k}}}}.
\end{equation}
\end{theorem}
\begin{proof}
To obtain (\ref{eqn:IIk21gf}), we count tilings of a $1 \times n$ rectangle with tiles of size $1\times 1,\ldots, 1 \times (k-1)$ according to the length of the right-most tile.
The generating function for the empty tiling is $1$.
The generating function for those tilings whose right-most tile has length $j \ge 1$ is $x^j y^{{{j}\choose{k}}} \II_{k\ldots 21}(x,y)$.
Combining these observations, we find that
$$\II_{k\ldots 21}(x,y) = 1 + \sum_{j\ge 1} x^j y^{{{j}\choose{k}}} \II_{k\ldots 21}(x,y).$$
Solve this equation for $\II_{k\ldots 21}(x,y)$ to obtain (\ref{eqn:IIk21gf}).
\end{proof}

For $r \le k$ one can now obtain $\II_{k\ldots 21}^r(x)$ from (\ref{eqn:IIk21gf}) by expanding the right side in powers of $y$ and finding the coefficient of $y^r$.
In lieu of this calculation, we use a combinatorial approach to obtain $\II_{k\ldots 21}^r(x)$.

\begin{theorem}
\label{thm:Uk21x}
For all $k$ and $r$ such that $0\leq r\leq k$ and all $n \ge 0$ we have
\begin{equation}
\label{eqn:Inrk21}
|I_n^r(231| k\ldots 21)| = \sum_{s_0,\ldots,s_r} \prod_{i=0}^r F_{k-1,s_i+1},
\end{equation}
where the sum on the right is over all sequences $s_0, \ldots, s_r$ of nonnegative integers such that $\sum_{i=0}^r s_i = n-kr.$
Moreover,
\begin{equation}
\label{eqn:Inrk21gf}
\II_{k\dots 21}^r(x)=\frac{x^{kr}}{(1-x-\dots-x^{k-1})^{r+1}}.
\end{equation}
\end{theorem}
\begin{proof}
Since $r \le k$, the tilings which correspond to the permutations in $I_n^r(231| k\ldots 21)$ are exactly those tilings of a $1\times n$ rectangle which contain precisely $r$ tiles of size $1 \times k$ and no tiles of length $k+1$ or more.
To build such a tiling, first order the tiles of size $1 \times k$;  there is one way to do this.
Now fix the sizes $s_0,\ldots,s_r$ of the gaps between these tiles;  observe that we must have $\sum_{i=0}^r s_i = n-kr$.
Finally, tile each of these gaps with tiles of length at most $k-1$;  there are $\prod_{i=0}^r F_{k-1,s_i+1}$ ways to do this.
Combine these observations to obtain (\ref{eqn:Inrk21}).

The proof of (\ref{eqn:Inrk21gf}) is similar to the proof of (\ref{eqn:Inrk21}).
\end{proof}

Using the same combinatorial techniques, we now compute $\II^1_\tau(x)$ for any layered permutation $\tau$.

\begin{theorem}
\label{thm:U1xlayered}
Fix a layered permutation $[l_1,\ldots,l_m]$.
Set $k_0 = l_1 - 1$, $k_m = l_m - 1$, and $k_i = \min(l_i-1, l_{i+1}-1)$ for $1 \le i \le m-1$.
Then we have
\begin{equation}
\label{eqn:I1layered}
|I_n^1(231| [l_1,\ldots,l_m])| = \sum_{s_0,\ldots,s_m} \prod_{i=0}^m F_{k_i,s_i+1},
\end{equation}
where the sum is over all sequences of nonnegative integers such that $\sum_{i=0}^ms_i = n-\sum_{i=1}^m l_i$.
Moreover,
\begin{equation}
\label{eqn:I1layeredgf}
\II^1_{[l_1,\ldots,l_m]}(x) = x^{\sum_{i=1}^m l_i} \prod_{i=0}^m \frac{1}{1-x-x^2-\ldots-x^{k_i}}.
\end{equation}
\end{theorem}
\begin{proof}
This is similar to the proof of Theorem \ref{thm:Uk21x}.
\end{proof}

Using Theorem \ref{thm:U1xlayered}, we now highlight a family of involutions which are enumerated by $k$-generalized Fibonacci numbers.

\begin{corollary}
\label{cor:kgenenum}
For all $k \ge 1$, all $l \ge 1$, and all $n \ge 0$, the number of involutions in $I_n(231)$ which contain exactly one subsequence of type $[1^k, l]$ is given by $F_{l-1,n-k-l+1}$.
\end{corollary}
\begin{proof}
Set $m = k+1$, $l_i = 1$ for $1 \le i \le k$, and $l_{k+1} = l$ in Theorem \ref{thm:U1xlayered} and observe that $F_{0,1} = 1$ and $F_{0,n} = 0$ whenever $n \neq 1$.
\end{proof}

\section{Involutions Which Contain 231 Once and Avoid Another Pattern}

In this section we consider those involutions in $I_n$ which contain exactly one subsequence of type 231 and avoid an additional pattern.
We begin by setting some notation.

\begin{definition}
For all $n \ge 0$ and any permutation $\tau$, we write $C_n(\tau)$ to denote the set of involutions in $I_n$ which avoid $\tau$ and which contain exactly one subsequence of type 231.
We write $\JJ_\tau(x)$ to denote the generating function given by
\begin{displaymath}
\JJ_\tau(x) = \sum_{n=0}^\infty |C_n(\tau)| x^n.
\end{displaymath}
\end{definition}

As we did for $I_n(231)$ in Section \ref{sec:In231tau}, we give a constructive bijection between the set of involutions of length $n$ which contain exactly one subsequence of type 231 and a certain set of tilings of a rectangle of size $1 \times n$.

\begin{proposition}
\label{prop:redblue}
Fix $n \ge 0$.
Then there exists a constructive bijection between the set of involutions in $I_n$ which contain exactly one subsequence of type 231 and the set of tilings of a $1 \times n$ rectangle using exactly one red rectangle of size $1 \times 4$ and blue rectangles of size $1 \times 1$, $1 \times 2, \ldots$.
\end{proposition}
\begin{proof}
Suppose we are given such a tiling;  we construct the corresponding involution as follows.
First number the squares $1,2,\ldots,n$ from left to right.
In each blue tile, reverse the order of the entries.
If the left-most entry of the red tile is $a$, then put the entries of the red tile in the order $a+3,a+1,a+2,a$.
It is routine to verify that the resulting permutation is an involution which contains exactly one subsequence of type 231 and that the given map is injective.
Therefore it is sufficient to show that the given map is surjective.
To do this, fix $\pi \in I_n$;  we argue by induction on $n$.
It is routine to verify the result when $n \le 4$, so we assume $n \ge 5$ and that the result holds for all $k \le n-1$.
Fix $j$ such that $\pi(j) = n$;  since $\pi$ is an involution we also have $\pi(n) = j$.
We consider two cases:  either $n$ appears in the subsequence of type 231 or $n$ does not appear in the subsequence of type 231.

If $n$ does not appear in the subsequence of type 231 then there exist permutations $\pi_1$ and $\pi_2$, exactly one of which contains exactly one subsequence of type 231, such that $\pi = \pi_1,n,\tilde{\pi_2}$.
Here $\tilde{\pi_2}$ is the sequence obtained by adding $j-1$ to every entry of $\pi_2$.
In this case the result follows from our bijection involving $I_n(231)$ and induction.

Now suppose the subsequence given by $\pi(a),n,\pi(b)$ has type 231.
If $\pi(b)>j$ then $\pi(a)>\pi(b)>j$ and the subsequence $\pi(a),n,j$ of $\pi$ is a second subsequence of type 231, which is a contradiction.
If $\pi(b) = j$ then $\pi(a)>j$.
Therefore, since $\pi$ is an involution, $a$ appears to the right of $n$ and $a < \pi(a)$.
It follows that $\pi(a),n,a$ is a second subsequence of type 231, which is a contradiction.
If $\pi(b)<j$ then $b$ appears to the left of $n$, since $\pi$ is an involution.
It follows that $b,n,j$ is a second subsequence of type 231, which is a contradiction.
Combining these observations, we find that $n$ does not appear in the subsequence of type 231.

It follows that the given map is surjective, as desired.
\end{proof}

Using Proposition \ref{prop:redblue}, we now enumerate those involutions in $I_n$ which contain exactly one subsequence of type 231.

\begin{theorem}
\label{thm:CnE}
The number of the involutions in $S_n$ which contain exactly one subsequence of type $231$ is given by $(n-1) 2^{n-6}$ for all $n\geq 5$.
\end{theorem}
\begin{proof}
Count the corresponding tilings in three cases:  the red tile is at the far left, the red tile is at the far right, or there are blue tiles on both sides of the red tile.
There are $2^{n-5}$ tilings of the first type, $2^{n-5}$ tilings of the second type, and $\sum_{i=1}^{n-5} 2^{i-1} 2^{n-i-5} = (n-5) 2^{n-6}$ tilings of the third type.
Combine these observations to obtain the desired result.
\end{proof}

We now use Proposition \ref{prop:redblue} to enumerate $C_n(k\ldots 21)$.
We observe that every involution in $S_n$ which contains exactly one pattern of type 231 contains a pattern of type 4231, and therefore a pattern of type 321.
With this in mind, we consider $C_n(k\ldots 21)$ for $k \ge 4$.

\begin{theorem}
\label{thm:Ck21}
Fix $k \ge 4$.
Then for all $n \ge 0$ we have
\begin{equation}
\label{eqn:Ck21enum}
|C_n(k\ldots 21)| = \sum_{i=0}^{n-4} F_{k-1,i+1} F_{k-1,n-i-3}.
\end{equation}
Moreover,
\begin{equation}
\label{eqn:Ck21gf}
\JJ_{k\ldots 21}(x) = \frac{x^4}{(1-x-\dots-x^{k-1})^2}.
\end{equation}
\end{theorem}
\begin{proof}
This is similar to the proof of Theorem \ref{thm:Uk21x}, using Proposition \ref{prop:redblue}.
\end{proof}

We conclude this section by giving a recursive procedure for computing $\JJ_\tau(x)$ when $\tau$ is a layered permutation.

\begin{theorem}
\label{thm:Clayered}
Fix positive integers $l_1,\ldots,l_m$ and $l$.
If $l \le 3$ then we have
\begin{equation}
\label{eqn:Clayeredl3}
\JJ_{[l_1,\ldots,l_m,l]}(x) = \frac{1}{1-x-\ldots-x^{l-1}}\left(\frac{x^l}{1-x}\JJ_{[l_1,\ldots,l_m]}(x)+x^4\II_{[l_1,\ldots,l_m]}(x)\right).
\end{equation}
If $l \ge 4$ then we have
\begin{equation}
\label{eqn:Clayeredl4}
\JJ_{[l_1,\ldots,l_m,l]}(x) = \frac{1}{1-x-\ldots-x^{l-1}}\left(\frac{x^l}{1-x}\JJ_{[l_1,\ldots,l_m]}(x)+x^4\II_{[l_1,\ldots,l_m,l]}(x)\right).
\end{equation}
\end{theorem}
\begin{proof}
This is similar to the proof of Theorem \ref{thm:Uk21x}, using Proposition \ref{prop:redblue}.
\end{proof}

We remark that it is clear from results given in \cite{Mansour33k} and (\ref{eqn:Ck21gf}) -- (\ref{eqn:Clayeredl4}) that for any layered permutation $\tau$, the set $C_n(\tau)$ can be enumerated in terms of $k$-generalized Fibonacci numbers.

\section{Involutions Which Contain 231 Once and Contain Another Pattern}

In this section we consider those involutions in $I_n$ which contain exactly one subsequence of type 231 and which also contain a given pattern.
We begin by setting some notation.

\begin{definition}
For all $n \ge 0$, all $r \ge 0$, and any permutation $\tau$, we write $C_n^r(\tau)$ to denote the set of involutions in $I_n$ which contain exactly one subsequence of type 231 and exactly $r$ subsequences of type $\tau$.
We write $\JJ_\tau^r(x)$ to denote the generating function given by
\begin{displaymath}
\JJ_\tau^r(x) = \sum_{n=0}^\infty |C_n^r(\tau)| x^n.
\end{displaymath}
We write $\JJ_\tau(x,y)$ to denote the generating function given by
\begin{displaymath}
\JJ_\tau(x,y) = \sum_{n=0}^\infty \sum_{r=0}^\infty |C_n^r(\tau)| x^n y^r.
\end{displaymath}
\end{definition}

As with $\II_\tau(x,y)$, finding a closed form expression for $\JJ_\tau(x,y)$ is difficult for general $\tau$, so we content ourselves here with a closed form expression for $\JJ_{k\ldots 21}(x,y)$.

\begin{theorem}
\label{the1}
For all $k\geq 4$ we have
\begin{equation}
\label{eqn:JJxy}
\JJ_{k\dots 21}(x,y)=\frac{x^4}{\left(1-\sum_{j\geq 1} x^jy^{{{j}\choose{k}}}\right)^2}.
\end{equation}
\end{theorem}
\begin{proof}
This is similar to the proof of Theorem \ref{thm:avoid231contain12kgf}, using Proposition \ref{prop:redblue}.
\end{proof}

For $r \le k$ one can now obtain $\JJ_{k\ldots 21}^r(x)$ from (\ref{eqn:JJxy}) by expanding the right side in powers of $y$ and finding the coefficient of $y^r$.
In lieu of this calculation, we use a combinatorial approach to obtain $\JJ_{k\ldots 21}^r(x)$.

\begin{theorem}
For all $k$ and $r$ such that $0 \le r \le k$ and all $n \ge 0$ we have
\begin{equation}
\label{eqn:Cnrk21}
|C_n^r(k\ldots 21)| = (r+1) \sum_{s_0,\ldots,s_{r+1}} \prod_{i=0}^{r+1} F_{k-1,s_i+1},
\end{equation}
where the sum on the right is over all sequences $s_0, \ldots, s_{r+1}$ of nonnegative integers such that $\sum_{i=0}^{r+1} s_i = n-kr-4$.
Moreover,
\begin{equation}
\label{eqn:JJk21gf}
\JJ_{k\ldots 21}^r(x) = \frac{(r+1) x^{kr+4}}{(1-x-x^2-\ldots-x^{k-1})^{r+2}}.
\end{equation}
\end{theorem}
\begin{proof}
This is similar to the proof of Theorem \ref{thm:Uk21x}, using Proposition \ref{prop:redblue}.
\end{proof}

We conclude the paper by giving a recursive procedure for computing $\JJ_\tau^1(x)$ for any layered permutation $\tau$.

\begin{theorem}
Fix positive integers $l_1,\ldots,l_m$ and $l$.
If $l \le 3$ then we have
\begin{equation}
\label{eqn:Crlayeredl3}
\JJ^1_{[l_1,\ldots,l_m,l]}(x) = \frac{x^l}{1-x-\dots-x^{l-1}}\JJ_{[l_1,\ldots,l_m]}^1(x).
\end{equation}
If $l \ge 4$ then we have
\begin{equation}
\label{eqn:Crlayeredl4}
\JJ^1_{[l_1,\ldots,l_m,l]}(x) = \frac{1}{1-x-\dots-x^{l-1}}\left(x^l\JJ_{[l_1,\dots,l_m]}^1(x)+x^4\II_{[l_1,\dots,l_m,l]}^1(x)\right).
\end{equation}
\end{theorem}
\begin{proof}
This is similar to the proof of Theorem \ref{thm:avoid231contain12kgf}, using Proposition \ref{prop:redblue}.
\end{proof}

We remark that it is clear from (\ref{eqn:I1layeredgf}) and (\ref{eqn:JJk21gf}) -- (\ref{eqn:Crlayeredl4})  that for any layered permutation $\tau$, the set $C^r_n(\tau)$ can be enumerated in terms of $k$-generalized Fibonacci numbers.

\end{document}